\documentclass[11pt,twoside]{article} 

\date{} 

\begin{document} 

%\centerline{\bf Journal's Title, Vol. x, 2018, no. xx, xxx - xxx} 

\centerline{} 

\centerline{} 

\centerline {\Large{\bf The asymptotics of the generalised}}

\centerline{}

\centerline {\Large{\bf Bessel function}}

\centerline{} 

\centerline{\bf {R. B. Paris}} 

\centerline{} 

\centerline{Division of Computing and Mathematics,} 

\centerline{Abertay University, Dundee DD1 1HG, UK} 

%\centerline{Address of Author2 third line} 

%\centerline{Address of Author2 fourth line} 
%
%\begin{document}
\def\f#1#2{\mbox{${\textstyle \frac{#1}{#2}}$}}
\def\dfrac#1#2{\displaystyle{\frac{#1}{#2}}}
\def\boldal{\mbox{\boldmath $\alpha$}}
{\newcommand{\Sgoth}{S\;\!\!\!\!\!/}
\newcommand{\bee}{\begin{equation}}
\newcommand{\ee}{\end{equation}}
\newcommand{\la}{\lambda}
\newcommand{\ka}{\kappa}
\newcommand{\al}{\alpha}
\newcommand{\fr}{\frac{1}{2}}
\newcommand{\fs}{\f{1}{2}}
\newcommand{\g}{\Gamma}
\newcommand{\br}{\biggr}
\newcommand{\bl}{\biggl}
\newcommand{\ra}{\rightarrow}
\newcommand{\gl}{\raisebox{-.8ex}{\mbox{$\stackrel{\textstyle >}{<}$}}}
\newcommand{\gtwid}{\raisebox{-.8ex}{\mbox{$\stackrel{\textstyle >}{\sim}$}}}
\newcommand{\ltwid}{\raisebox{-.8ex}{\mbox{$\stackrel{\textstyle <}{\sim}$}}}
\renewcommand{\topfraction}{0.9}
\renewcommand{\bottomfraction}{0.9}
\renewcommand{\textfraction}{0.05}
\newcommand{\mcol}{\multicolumn}
%\date{}
%\maketitle
%\pagestyle{myheadings}
%\markboth{\hfill \sc R. B.\ Paris  \hfill}
%{\hfill \sc  The generalised Bessel function \hfill}
\begin{abstract}
We demonstrate how the asymptotics for large $|z|$ of the generalised Bessel function 
\[{}_0\Psi_1(z)=\sum_{n=0}^\infty\frac{z^n}{\Gamma(an+b) n!},\]
where $a>-1$ and $b$ is any number (real or complex), may be obtained by exploiting the well-established asymptotic theory of the generalised Wright function ${}_p\Psi_q(z)$. A summary of this theory is given and an algorithm for determining the coefficients in the associated exponential expansions is discussed in an appendix. We pay particular attention to the case $a=-\fs$, where the expansion for $z\to\pm\infty$ consists of an exponentially small contribution that undergoes a Stokes phenomenon. 

We also examine the different nature of the asymptotic expansions as a function of $\arg\,z$ when $-1<a<0$, taking into account the Stokes phenomenon that occurs on the rays $\arg\,z=0$ and $\arg\,z=\pm\pi(1+a)$ for the associated function ${}_1\Psi_0(z)$. These regions are more precise than those given by Wright in his 1940 paper. Numerical computations are carried out to verify several of the expansions developed in the paper.
\vspace{0.4cm}

\vspace{0.3cm}

\end{abstract}

\vspace{0.3cm}

%\noindent $\,$\hrulefill $\,$}

%\begin{center}
{\bf 1. \  Introduction}
%\end{center}
\setcounter{section}{1}
\setcounter{equation}{0}
\renewcommand{\theequation}{\arabic{section}.\arabic{equation}}
The generalised Bessel function is defined by
\bee\label{e10}
{}_0\Psi_1(z)\equiv {}_0\Psi_1\bl(\begin{array}{c}-\!-\\(a,b)\end{array}\!\!\left|\,z\br)\right.=\sum_{n=0}^\infty\frac{z^n}{\g(an+b) n!},
\ee
where $a$ is supposed real and $b$ is an arbitrary complex parameter. The series 
converges for all finite $z$ provided $a>-1$ and, when $a=1$, it reduces to the modified Bessel function $z^{(1-b)/2}I_{b-1}(2\sqrt{z})$. The asymptotics of this function were first studied by Wright \cite{W34, W39} using the method of steepest descents applied to a suitable integral representation. A more recent investigation by Wong and Zhao \cite{WZ1, WZ2} considered the hyperasymptotic expansion of ${}_0\Psi_1(z)$ and the associated Stokes phenomenon. The case $b=0$ finds application in probability theory and is discussed in \cite{PV3}, where it is referred to as a `reduced' Wright function.

The series in (\ref{e10}) is a particular case of the generalised Wright function ${}_p\Psi_q(z)$ defined below in (\ref{e20}). A brief presentation of the large-$z$ asymptotic theory of this function is given in Section 2.
In this paper we derive and summarise the asymptotic expansion of the generalised Bessel function using the asymptotic theory of the generalised Wright function. It is found that 
the asymptotic expansion of ${}_0\Psi_1(z)$ as $|z|\to\infty$ separates into two distinct cases according as $a>0$ and $-1<a<0$. Our principal aim is to reconsider the exponentially small contribution to the expansion in the particular case $a=-\fs$ as $z\to\pm\infty$. This was discussed by Wright \cite[\S 5]{W39} though, as we shall show, his form of the exponentially small contribution is not correct as he took no account of the Stokes phenomenon.

We also examine the different nature of the asymptotic expansions as a function of $\arg\,z$ when $-1<a<0$, taking into account the Stokes phenomenon that occurs on the rays $\arg\,z=0$ and $\arg\,z=\pm\pi(1+a)$ for the associated function ${}_1\Psi_0(z)$ defined in (\ref{e36}) below . These regions are more precise than those given by Wright in his 1940 paper \cite{W39}; for a summary of these expansions, see also \cite{WZ2}.  

\vspace{0.6cm}

%\begin{center}
{\bf 2. \  The asymptotic expansion of ${}_p\Psi_q(z)$ for $|z|\to\infty$}
%\end{center}
\setcounter{section}{2}
\setcounter{equation}{0}
\renewcommand{\theequation}{\arabic{section}.\arabic{equation}}
The generalised Wright function is defined by the series
\bee\label{e20}
{}_p\Psi_q(z)\equiv{}_p\Psi_q\bl(\!\!\begin{array}{c}(\alpha_1,a_1), \ldots ,(\alpha_p,a_p)\\(\beta_1, b_1), \ldots ,(\beta_q,b_q)\end{array}\!\!\left|\,z\br)\right.=\sum_{n=0}^\infty g(n)\,\frac{z^n}{n!}, 
\ee
\bee\label{e20a}
g(n)=\frac{\prod_{r=1}^p\Gamma(\alpha_rn+a_r)}{\prod_{r=1}^q\Gamma(\beta_rn+b_r)},
\ee
where $p$ and $q$ are nonnegative integers, the parameters $\alpha_r$  and 
$\beta_r$ are real and positive and $a_r$ and $b_r$ are
arbitrary complex numbers. We also assume that the $\alpha_r$ and $a_r$ are subject to 
the restriction
\bee\label{e20ab}
\alpha_rn+a_r\neq 0, -1, -2, \ldots \qquad (n=0, 1, 2, \ldots\ ;\, 1\leq r \leq p)
\ee
so that no gamma function in the numerator in (\ref{e20}) is singular. The determination of the asymptotic expansion of ${}_p\Psi_q(z)$ for $|z|\ra\infty$ and finite 
values of the parameters has a long history. Detailed investigations were carried out by Wright \cite{W40} and by
Braaksma \cite{Br} for a more general class of integral functions than (\ref{e20}). We present below a summary of the main expansion theorems related to the asymptotics of ${}_p\Psi_q(z)$ for large $|z|$; for a recent presentation, see \cite{P17}.

We introduce the parameters associated\footnote{Empty sums and products are to be interpreted as zero and unity, respectively.} with $g(n)$ given by
\[\kappa=1+\sum_{r=1}^q\beta_r-\sum_{r=1}^p\alpha_r, \qquad 
h=\prod_{r=1}^p\alpha_r^{\alpha_r}\prod_{r=1}^q\beta_r^{-\beta_r},\]
\bee\label{e21}
\vartheta=\sum_{r=1}^pa_r-\sum_{r=1}^qb_r+\f{1}{2}(q-p),\qquad \vartheta'=1-\vartheta.
\ee
If it is supposed that $\alpha_r$ and $\beta_r$ are such that $\kappa>0$ then ${}_p\Psi_q(z)$ 
is uniformly and absolutely convergent for all finite $z$. If $\kappa=0$, the sum in (\ref{e20})
has a finite radius of convergence equal to $h^{-1}$, whereas for $\kappa<0$ the sum is divergent 
for all nonzero values of $z$. The parameter $\kappa$ will be found to play a critical role 
in the asymptotic theory of ${}_p\Psi_q(z)$ by determining the sectors in the $z$-plane 
in which its behaviour is either exponentially large, algebraic or exponentially small 
in character as $|z|\ra\infty$. 

We first introduce the exponential expansion $E_{p,q}(z)$ and the 
algebraic expansion $H_{p,q}(z)$ associated with ${}_p\Psi_q(z)$.
The exponential expansion is given by the formal asymptotic sum
\bee\label{e22c}
E_{p,q}(z):=Z^\vartheta e^Z\sum_{j=0}^\infty A_jZ^{-j}, \qquad Z=\kappa (hz)^{1/\kappa},
\ee
where the coefficients $A_j$ are those appearing in the inverse factorial expansion of $g(s)/s!$ given by  
\bee\label{e22a}
\frac{g(s)}{\g(1+s)}=\kappa (h\kappa^\kappa)^s\bl\{\sum_{j=0}^{M-1}\frac{A_j}{\Gamma(\kappa s+\vartheta'+j)}
+\frac{\rho_M(s)}{\Gamma(\kappa s+\vartheta'+M)}\br\}.
\ee
Here $g(s)$ is defined in (\ref{e20a}) with $n$ replaced by $s$, $M$ is a positive integer and $\rho_M(s)=O(1)$ for $|s|\ra\infty$ in $|\arg\,s|<\pi$.
The leading coefficient $A_0$ is specified by
\bee\label{e22b}
A_0=(2\pi)^{\frac{1}{2}(p-q)}\kappa^{-\frac{1}{2}-\vartheta}\prod_{r=1}^p
\alpha_r^{a_r-\frac{1}{2}}\prod_{r=1}^q\beta_r^{\frac{1}{2}-b_r}.
\ee
The coefficients $A_j$ are independent of $s$ and depend only on the parameters $p$, $q$, $\alpha_r$, 
$\beta_r$, $a_r$ and $b_r$. An algorithm for their evaluation is described in Appendix A.

The algebraic expansion $H_{p,q}(z)$ follows from the Mellin-Barnes integral representation \cite[ \S 2.4]{PK}
\bee\label{e24aa}
{}_p\Psi_q(z)=\frac{1}{2\pi i}\int_{-\infty i}^{\infty i} \Gamma(s)g(-s)(ze^{\mp\pi i})^{-s}ds,\qquad |\arg(-z)|<\pi(1-\fs\kappa),
\ee
where the path of integration is indented near $s=0$ to separate\footnote{This is always 
possible when the condition (\ref{e20ab}) is satisfied.} the poles of $\g(s)$ at $s=-k$ from those of 
$g(-s)$ situated at 
\bee\label{e24a}
s=(a_r+k)/\alpha_r, \qquad k=0, 1, 2, \dots\, \ (1\leq r\leq p).
\ee
In general there will be $p$ such sequences of simple poles though, depending on the values 
of $\alpha_r$ and $a_r$, some of these poles could be multiple poles or even ordinary 
points if any of the $\Gamma(\beta_rs+b_r)$ are singular there. Displacement of the contour to the 
right over the poles of $g(-s)$ then yields the algebraic expansion of 
${}_p\Psi_q(z)$ valid in the sector in (\ref{e24aa}). 

If it is assumed that the parameters are such that 
the poles in (\ref{e24a}) are all simple we obtain the algebraic expansion given by 
$H_{p,q}(z)$, where
\bee\label{e25}
H_{p,q}(z):=\sum_{m=1}^p\alpha_m^{-1}z^{-a_m/\alpha_m}S_{p,q}(z;m)
\ee
and $S_{p,q}(z;m)$ denotes the formal asymptotic sum
\bee\label{e25a}
S_{p,q}(z;m):=\sum_{k=0}^\infty \frac{(-)^k}{k!}\Gamma\left(\frac{k+a_m}{\alpha_m}\right)\,
\frac{\prod_{r=1}^{'\,p}\Gamma(a_r-\alpha_r(k+a_m)/\alpha_m)}
{\prod_{r=1}^q\Gamma(b_r-\beta_r(k+a_m)/\alpha_m)} z^{-k/\alpha_m},
\ee
with the prime indicating the omission of the term corresponding to $r=m$ in the product. 
This expression in (\ref{e25}) consists of (at most) $p$ expansions each with the leading behaviour 
$z^{-a_m/\alpha_m}$ ($1\leq m\leq p$).
When the parameters $\alpha_r$ and $a_r$ are such that some of the poles 
are of higher order, the expansion (\ref{e25a}) is invalid and the residues must 
then be evaluated according to the multiplicity of the poles concerned; this will lead to terms involving $\log\,z$ in the algebraic expansion.

The expansion theorems for ${}_p\Psi_q(z)$ are as follows. Throughout we let $\epsilon$ denote an arbitrarily 
small positive quantity.
\newtheorem{theorem}{Theorem}
\begin{theorem}$\!\!\!.$
When $0<\kappa\leq 2$, then 
\bee\label{e24}
{}_p\Psi_q(z)\sim\left\{\begin{array}{ll}
E_{p,q}(z)+H_{p,q}(ze^{\mp\pi i}) & (|\arg\,z|\leq \min\{\pi-\epsilon,\pi\kappa-\epsilon\}) \\
%\\ \fs {\cal E}_{p,q}(z)+H_{p,q}^o(ze^{\mp\pi i}) & \mbox{on} & \arg\,z=\pm\pi\kappa\\ 
\\H_{p,q}(ze^{\mp\pi i}) & 
 (\pi\kappa+\epsilon\leq |\arg\,z|\leq\pi,\ 0<\kappa<1)\\
\\E_{p,q}(z)+E_{p,q}(ze^{\mp2\pi i}) & \\
+H_{p,q}(ze^{\mp\pi i}) & 
 (|\arg\,z|\leq\pi,\ 1<\kappa\leq 2)\end{array} \right.
\ee
as $|z|\ra\infty$. The upper or lower signs 
are chosen according as $\arg\,z>0$ or $\arg\,z<0$, respectively. 
\end{theorem}
\begin{theorem}$\!\!\!.$
When $\kappa>2$ we have
%\footnote{In \cite{W1}, the expansion was given in terms of the two dominant expansions only, viz. $E_{p,q}(z)$ and $E_{p,q}(ze^{\mp2\pi i})$, corresponding to $n=0$ and $n=\pm 1$ in (\ref{e23}).} 
\bee\label{e23}
{}_p\Psi_q(z)\sim \sum_{n=-N}^N E_{p,q}(ze^{2\pi in})+H_{p,q}(ze^{\mp\pi iz})
\ee
as $|z|\to\infty$ in the sector $|\arg\,z|\leq\pi$. The integer $N$ is chosen such that it is the smallest integer satisfying $2N+1>\fs\kappa$ and the upper or lower is chosen according as $\arg\,z>0$ or $\arg\,z<0$, respectively.

In this case the asymptotic behaviour of ${}_p\Psi_q(z)$ is exponentially large for all values of $\arg\,z$ and, consequently, the algebraic expansion may be neglected. The sums $E_{p,q}(ze^{2\pi in})$ are exponentially large (or oscillatory) as $|z|\to\infty$ for values of $\arg\,z$ satisfying $|\arg\,z+2\pi n|\leq\fs\pi\kappa$. 
\end{theorem}

The division of the $z$-plane into regions where ${}_p\Psi_q(z)$ possesses exponentially large or algebraic behaviour for large $|z|$ is illustrated in Fig.~1.
When $0<\kappa<2$, the exponential expansion $E_{p,q}(z)$ is still present in the sectors $\fs\pi\kappa<|\arg\,z|<\min\{\pi,\pi\kappa\}$, where it is subdominant. The rays $\arg\,z=\pm\pi\kappa$ ($0<\kappa<1$), where $E_{p,q}(z)$ is {\it maximally} subdominant with respect to $H_{p,q}(ze^{\mp\pi i})$, are called Stokes lines.\footnote{The positive real axis $\arg\,z=0$ is also a Stokes line where the algebraic expansion is maximally subdominant.} As these rays are crossed (in the sense of increasing $|\arg\,z|$) the exponential expansion switches off according to Berry's now familiar error-function smoothing law \cite{B}; see \cite{P10} for details. The rays $\arg\,z=\pm\fs\pi\kappa$, where $E_{p,q}(z)$ is oscillatory and comparable to $H_{p,q}(ze^{\mp\pi i})$, are called anti-Stokes lines.
\begin{figure}[t]
\centering
\begin{picture}(200,200)(0,0)
\put(0,100){\line(1,0){200}}
\put(100,100){\line(-1,1){80}}
\put(100,100){\line(-1,-1){80}}
\multiput(100,100)(14,14){6}{\line(1,1){10}}
\multiput(100,100)(14,-14){6}{\line(1,-1){10}}
\put(125,110){\footnotesize{Exponentially large}}
\put(130,85){\footnotesize{$+$ Algebraic}}
\put(60,160){\footnotesize{Exponentially small}}
\put(70,140){\footnotesize{$+$ Algebraic}}
\put(60,40){\footnotesize{Exponentially small}}
\put(70,20){\footnotesize{$+$ Algebraic}}
\put(20,110){\footnotesize{Algebraic}}
\put(177,160){\footnotesize{$\theta=\pi\kappa/2$}}
\put(185,180){\footnotesize{anti-Stokes line}}
\put(177,40){\footnotesize{$\theta=-\pi\kappa/2$}}
\put(185,20){\footnotesize{anti-Stokes line}}
\put(-5,160){\footnotesize{$\theta=\pi\kappa$}}
\put(10,190){\footnotesize{Stokes line}}
\put(-5,40){\footnotesize{$\theta=-\pi\kappa$}}
\put(10,8){\footnotesize{Stokes line}}
\end{picture}
\caption{\small{The exponentially large and algebraic sectors associated with ${}_p\Psi_q(z)$ in the complex $z$-plane with $\theta=\arg\,z$ when $0<\kappa<1$. The Stokes and anti-Stokes lines are indicated.}}
\end{figure}
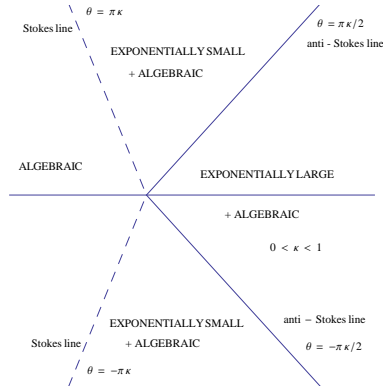

We omit the expansion {\it on} the Stokes lines $\arg\,z=\pm\pi\kappa$ in Theorem 1; the details in the case $p=1$, $q\geq0$ are discussed in \cite{P14}. We shall return to this point in Section 4 when discussing the case $a=-\fs$.
Since $E_{p,q}(z)$ is exponentially small in $\fs\pi\kappa<|\arg\,z|\leq\min\{\pi, \pi\kappa\}$, then in the sense of Poincar\'e, the expansion $E_{p,q}(z)$ can be neglected in these sectors.
Similarly, $E_{p,q}(ze^{-2\pi i})$ is exponentially small compared to $E_{p,q}(z)$ in $0\leq\arg\,z<\pi$ and so can be neglected when $1<\kappa<2$. However, in the vicinity of $\arg\,z=\pi$, these last two expansions are of comparable magnitude and, for real parameters, they combine to generate a real result on this ray. A similar remark applies to $E_{p,q}(ze^{2\pi i})$ in $-\pi<\arg\,z\leq 0$. 

\vspace{0.6cm}

%\begin{center}
{\bf 3. \  The asymptotic expansion of ${}_0\Psi_1(z)$ when $a>0$}
%\end{center}
\setcounter{section}{3}
\setcounter{equation}{0}
\renewcommand{\theequation}{\arabic{section}.\arabic{equation}}
The generalised Bessel function ${}_0\Psi_1(z)$ in (\ref{e10}) is, from (\ref{e21}), associated with the parameters 
\[\kappa=1+a,\qquad h=a^{-a},\qquad \vartheta=\fs-b.\]
It follows that when $a>0$ we have $\kappa>1$. Since $p=0$, the algebraic expansion $H_{0,1}(z)\equiv 0$ and the exponential expansion is 
\bee\label{e301}
E_{0,1}(z)=Z^\vartheta e^Z \sum_{j=0}^\infty A_jZ^{-j},\qquad Z=\kappa(hz)^{1/\kappa},
\ee
where, from (\ref{e22b}), $A_0=(a/\kappa)^\vartheta /\sqrt{2\pi\kappa}$. The normalised coefficients $c_j\equiv A_j/A_0$ can be obtained by the algorithm described in Appendix A, with $g(s)=1/\g(as+b)$, to yield the first few coefficients given by
\[c_0=1,\quad c_1=-\frac{1}{24a}\bl\{(2+a)(1+2a)-12b(1+a-b)\br\},\]
\[c_2=\frac{1}{1152a^2}\bl\{(2+a)(1+2a)(2-19a+2a^2)+24b(1+a)(2+7a-6a^2)\]
\[\hspace{6cm}-24b^2(4-5a-20a^2)-96b^3(1+5a)+144b^4\br\},\]
\[c_3=\frac{1}{414720a^3}\bl\{(2+a)(1+2a)(556+1628a-9093a^2+1628a^3+556a^4)\]
\[-180b(1+a)(12-172a-417a^2+516a^3-20a^4)-180b^2(76+392a-567a^2-1288a^3\]
\[+364a^4)
+1440b^3(8-63a-147a^2+112a^3)+10800b^4(2+7a-14a^2)\]
\bee\label{e31c}
-8640b^5(1-7a)-8640b^6\br\}.
\ee
The coefficients $c_j$ ($1\leq j\leq 10$) for the case $a=\fs$, $b=\f{5}{4}$ are presented in Table 1.
\begin{table}[th]
\caption{\footnotesize{The normalised coefficients $c_j=A_j/A_0$ for ${}_0\Psi_1(z)$ in (\ref{e10}) when $a=\fs$ and $b=\f{5}{4}$. }}
\begin{center}
\begin{tabular}{c|l||c|l}
\mcol{1}{c|}{$j$} & \mcol{1}{c||}{$c_j$} & \mcol{1}{c|}{$j$} & \mcol{1}{c}{$c_j$}\\
[.05cm]\hline
&&& \\[-0.2cm]
1 & $-\f{5}{48}$ & 2 & $-\f{455}{4608}$\\
&&& \\[-0.2cm]
3 & $-\f{85085}{663552}$& 4 & $-\f{24079055}{127401984}$\\
&&& \\[-0.2cm]
5 & $-\f{1511535025}{6115295232}$ & 6 & $+\f{26957055125}{1761205026816}$\\
&&& \\[-0.2cm]
7 & $+\f{215144256952625}{84537841287168}$ & 8 & $+\f{570314645402376875}{32462531054272512}$\\
&&& \\[-0.2cm]
9 & $+\f{1304836837479714163625}{14023813415445725184}$ & 10 & $+\f{560395062780446967448375}{1346286087882789617664}$\\
[.15cm]\hline
\end{tabular}
\end{center}
\end{table}

Then from Theorems 1 and 2 we have 
\begin{theorem}$\!\!\!.$\ When $a>0$, we have the expansions of the generalised Bessel function
\bee\label{e32}
{}_0\Psi_1(z)\sim \left\{\begin{array}{ll} E_{0,1}(z)+E_{0,1}(ze^{\mp 2\pi i}) & (0<a\leq 1)\\
\\
\displaystyle{\sum_{n=-N}^N }E_{0,1}(ze^{2\pi in}) & (a>1)\end{array}\right.
\ee
for $|z|\to\infty$ in $|\arg\,z|\leq\pi$, where the upper or lower sign is chosen according as $\arg\,z>0$ or $\arg\,z<0$, respectively. The exponential expansion is defined in (\ref{e301}), with $A_0=(a/\kappa)^\vartheta/\sqrt{2\pi\kappa}$. The integer $N$ is the smallest integer satisfying $2N+1>\fs(1+a)$. 
\end{theorem}
We remark that when $a=1$ ($\kappa=2$) the expansion of ${}_0\Psi_1(z)$ can be obtained from that of the modified Bessel function.

From the discussion after Theorem 2, the expansions $E_{0,1}(ze^{\mp 2\pi i})$ are subdominant in the upper (lower) half-plane and switch on (in the sense of increasing $|\arg\,z|$) across the Stokes lines for these two functions; see \cite[\S 3.4]{P17} for a numerical example.
Since the exponential factors associated with $E_{0,1}(z)$ and $E_{0,1}(ze^{\mp2\pi i})$ are $\exp [|Z| e^{i\theta/\kappa})]$ and $\exp [|Z|e^{(\theta\mp2\pi i)/\kappa})]$, where $\theta=\arg\,z$, the greatest difference between the real parts of these factors occurs when
\[\sin \bl(\frac{\theta}{\kappa}\br)=\sin \bl(\frac{\theta\mp2\pi}{\kappa}\br);\]
that is, on the Stokes lines $\theta=\pm\pi(1-\fs\kappa)$. 
Consequently, when $1<\kappa<2$ ($0<a<1$), the expansions $E_{0,1}(ze^{\mp2\pi i})$ are not present on the positive real axis, but make a significant contribution on the negative real axis. Thus, on the real axis we have the expansions
\bee\label{e33}
{}_0\Psi_1(x) \sim E_{0,1}(x) = X^\vartheta e^X \sum_{j=0}^\infty A_jX^{-j} \qquad (0<a<1)
\ee
and
\[{}_0\Psi_1(-x) \sim E_{0,1}(xe^{\pi i})+E_{0,1}(xe^{-\pi i})\]
\bee\label{e34}
=2X^\vartheta e^{X\cos \pi/\kappa} \sum_{j=0}^\infty A_j X^{-j} \cos\,\bl[X\sin \frac{\pi}{\kappa}+\frac{\pi}{\kappa}(\vartheta-j)\br]\qquad (0<a\leq 1)
\ee
as $x\to+\infty$, where $X=\kappa(hx)^{1/\kappa}$.

An example when the parameter $a>1$ is furnished by the following. When $a=3$ ($\kappa=4$), we have $N=1$ and from the second expansion in Theorem 3
\begin{eqnarray*}
{}_0\Psi_1(x)&\sim& E_{0,1}(x)+E_{0,1}(xe^{2\pi i})+E_{0,1}(xe^{-2\pi i})\\
&=&X^\vartheta e^X \sum_{j=0}^\infty A_jX^{-j}+2X^\vartheta\sum_{j=0}^\infty A_jX^{-j} \cos\, \bl[X+\frac{\pi}{4}(2\vartheta-j)\br]
\end{eqnarray*}
as $x\to+\infty$. For negative argument, we find upon neglecting the exponentially small contribution $E_{0,1}(xe^{3\pi i})$,
\begin{eqnarray*}
{}_0\Psi_1(-x)&\sim& E_{0,1}(xe^{\pi i})+E_{0,1}(xe^{-\pi i})\\
&=& 2X^\vartheta e^{X\cos \pi/4} \sum_{j=0}^\infty A_jX^{-j} \cos\,\bl[X\sin\frac{\pi}{4}+\frac{\pi}{4}(\vartheta-j)\br].
\end{eqnarray*}
as $x\to+\infty$.
\vspace{0.6cm}

%\begin{center}
{\bf 4. \  The asymptotic expansion of ${}_0\Psi_1(z)$ when $-1<a<0$}
%\end{center}
\setcounter{section}{4}
\setcounter{equation}{0}
\renewcommand{\theequation}{\arabic{section}.\arabic{equation}}
Let $a=-\sigma$, with $0<\sigma<1$. Use of the reflection formula for the gamma function with $\vartheta=\fs-b$ yields
\begin{eqnarray}
{}_0\Psi_1(z)&=&\frac{1}{\pi}\sum_{n=0}^\infty \frac{\g(\sigma n\!+\!1\!-\!b)}{n!}\,z^n \sin \pi(-\sigma n+b)\nonumber\\
&=&\frac{1}{2\pi}\bl\{e^{\pi i\vartheta} {}_1\Psi_0(ze^{\pi i\sigma})+e^{-\pi i\vartheta} {}_1\Psi_0(ze^{-\pi i\sigma})\br\},\label{e35}
\end{eqnarray}
where ${}_1\Psi_0(z)\equiv {}_1\Psi_0((\sigma, 1-b);-\!-;z)$.
The expansion of ${}_0\Psi_1(z)$ with $-1<a<0$ can then be constructed from knowledge of the expansion of the associated function ${}_1\Psi_0(ze^{\pm\pi i\sigma})$. 
To keep the presentation as clear as possible, we restrict our attention in this section to the most commonly occurring  case of real $z$ in (\ref{e35}). We discuss the case of complex $z$ in Section 5. 
\vspace{0.4cm}

\noindent{\bf 4.1.\ The expansion of ${}_1\Psi_0(z)$}
\vspace{0.3cm}

\noindent
In this section we deal with the expansion of the associated function
\bee\label{e36}
{}_1\Psi_0(z)=\sum_{n=0}^\infty \frac{\g(\sigma n+\delta)}{n!}\,z^n \qquad (0<\sigma<1)
\ee
provided $\g(\sigma n+\delta)$ is regular for $n=0, 1, 2, \ldots\ $. The parameters associated with (\ref{e36}) are $\kappa=1-\sigma$, $h=\sigma^\sigma$ and $\vartheta=\delta-\fs$. From (\ref{e22c}) and (\ref{e25}), the algebraic and exponential expansions are
\bee\label{e37}
H_{1,0}(z)=\frac{1}{\sigma}\sum_{k=0}^\infty \frac{(-)^k}{k!} \g\bl(\frac{k+\delta}{\sigma}\br) z^{-(k+\delta)/\sigma}, \qquad
E_{1,0}(z)=Z^\vartheta e^Z\sum_{j=0}^\infty A_j(\sigma)Z^{-j},
\ee
where $Z$ defined in (\ref{e22c}). We shall find it convenient in the case of negative $a$ to denote the dependence of the coefficients on the parameter $\sigma$ in the exponential expansion  by writing $A_j(\sigma)$; we omit to denote the dependence on the parameter $\delta$. From (\ref{e22b}) the leading coefficient is
\bee\label{e37a}
A_0(\sigma)=(2\pi/\kappa)^{1/2}\,(\sigma/\kappa)^\vartheta.
\ee
Then, since $0<\kappa<1$, we obtain from Theorem 1 the large-$z$ expansion
\bee\label{e4exp}
{}_1\Psi_0(z)\sim \left\{\begin{array}{ll}E_{1,0}(z)+H_{1,0}(ze^{\mp\pi i}) & (|\arg\,z|\leq\pi\kappa-\epsilon) \\
\\
H_{1,0}(ze^{\mp\pi i}) & (\pi\kappa+\epsilon\leq|\arg\,z|\leq\pi),\end{array}\right.
\ee
where the upper or lower signs are chosen according as $\arg\,z>0$ or $\arg\,z<0$, respectively. With $g(s)=\g(\sigma s+\delta)$ for the function in (\ref{e36}), the algorithm in Appendix A shows that the first few normalised coefficients $c_j=A_j(\sigma)/A_0(\sigma)$ are\footnote{It can be shown that if the coefficients in (\ref{e31c}) are denoted by $c_j\equiv c_j(a,b)$, then the coefficients in (\ref{ecj}) when $\delta=1-b$ are given by $c_j(-\sigma,b)$.}
\[c_0=1\qquad c_1=\frac{1}{24\sigma}\{2+7\sigma+2\sigma^2-12\delta(1+\sigma)+12\delta^2\},\]
\[c_2=\frac{1}{1152\sigma^2}\{4+172\sigma+417\sigma^2+172\sigma^3+4\sigma^4-24\delta(6+41\sigma+41\sigma^2+6\sigma^3)\]
\[\hspace{6cm}+120\delta^2(4+11\sigma+4\sigma^2)-480\delta^3(1+\sigma)+144\delta^4\},\]
\[c_3=\frac{1}{414720 \sigma^3}\{(-1112 + 9636 \sigma + 163734 \sigma^2 + 336347 \sigma^3 + 
  163734 \sigma^4 + 9636 \sigma^5\]
  \[ - 1112 \sigma^6)-\delta(3600 + 220320 \sigma + 
  929700 \sigma^2 + 929700 \sigma^3  + 220320 \sigma^4  + 3600 \sigma^5)\] 
  \[+ \delta^2(65520  + 715680 \sigma + 1440180 \sigma^2 + 715680 \sigma^3 + 
  65520 \sigma^4)\] 
  \[ - \delta^3(161280  + 816480 \sigma  + 816480 \sigma^2  +
  161280 \sigma^3)\]
  \bee\label{ecj} +\delta^4 (151200  + 378000 \sigma + 151200 \sigma^2) - 
  60480\delta^5(1 + \sigma) + 8640 \delta^6\}.
  \ee
  
The expansion of ${}_1\Psi_0(z)$ on the Stokes lines $\arg\,z=\pm\pi\kappa$ is of a more recondite nature. This has been considered for the more general  
function ${}_1\Psi_q(xe^{\pm\pi i\kappa})$ for integer $q\geq 0$ in \cite[\S 5]{P14}. From (4.4), (4.19)--(4.24) of this reference with $q=0$ and $z=xe^{\pm\pi i\kappa}$, $x>0$ it is found that\footnote{There is a misprint in \cite[(4.24)]{P14}: the sign of $i$ should be reversed in both instances.}
\[{}_1\Psi_0(xe^{\pm\pi i\kappa})=\frac{e^{\pm\pi i\delta}}{\sigma} \sum_{k=0}^{m_o-1}\frac{\g(\frac{k+\delta}{\sigma})}{k!}\,x^{-(k+\delta)/\sigma}\]
\bee\label{e46}
+(Xe^{\pm\pi i})^\vartheta e^{-X} \bl\{\sum_{j=0}^{M-1}\bl(\frac{1}{2}A_j(\sigma)\pm \frac{i B_j(\sigma)}{\sqrt{2\pi X}}\br)(-X)^{-j}+O(X^{-M})\br\}
\ee
as $x\to+\infty$, where $X=\kappa(hx)^{1/\kappa}$ and $m_o$ denotes the optimal truncation index for the algebraic expansion. Some routine algebra shows that the optimal truncation index $m_o$ satisfies
\bee\label{e46a}
m_0=\frac{\sigma}{\kappa}(X+\f{3}{2})-\frac{1+2\delta}{2\kappa}+\alpha,
\ee
where $|\alpha|$ is bounded (and chosen such that $m_o$ is an integer).

The coefficients $A_j(\sigma)$ are those appearing in (\ref{e37}) and the $B_j(\sigma)$ are defined by
\bee\label{e47}
B_j(\sigma)=\sum_{k=0}^j (-2)^k (\fs)_k A_{j-k}(\sigma)\,G_{2k,j-k}(\sigma).
\ee
The coefficients $G_{k,j}(\sigma)$ appear in the expansion 
\bee\label{e49}
\frac{\mu\tau^{\gamma_j-1}}{1-\tau^\mu}\,\frac{d\tau}{dw}=-\frac{1}{w}+\sum_{k=0}^\infty G_{k,j}(\sigma)w^k,\qquad \fs w^2=\tau-\log\,\tau-1, 
\ee
where
\[\mu:=\frac{\kappa}{\sigma}=\frac{1-\sigma}{\sigma},\qquad \gamma_j:= \delta(1+\mu)-\frac{1}{2}+\mu m_o-j-X.\]
The branch of $w(\tau)$ is chosen such that $w\sim \tau-1$ as $\tau\ra 1$, so that upon reversion of the $w$-$\tau$ mapping we find
\bee\label{e49a}
\tau=1+w+\f{1}{3}w^2+\f{1}{36}w^3-\f{1}{270}w^4+\f{1}{4320}w^5+ \cdots\ .
\ee

For $0\leq k\leq 2$, we find from (\ref{e49}) and (\ref{e49a}) with the aid of {\it Mathematica}
\[G_{0,j}(\sigma)=-\gamma_j+\frac{1}{6}+\frac{1}{2}\mu,\]
\[G_{2,j}(\sigma)=-\frac{1}{6}\gamma_j^3+\frac{1}{4}(1+\mu)\gamma_j^2-\frac{1}{12}(1+3\mu+\mu^2)\gamma_j+\frac{1}{1080}(2+45\mu+45\mu^2),\]
\[G_{4,j}(\sigma)=-\frac{1}{120}\gamma_j^5+\frac{1}{144}(5+3\mu)\gamma_j^4-\frac{1}{216}(10+15\mu+3\mu^2)\gamma_j^3+\frac{1}{144}(3+10\mu+5\mu^2)\gamma_j^2\]
\bee\label{e42e}
-\frac{1}{4320}(5+90\mu+100\mu^2-6\mu^4)\gamma_j+\frac{1}{36288}(-13+21\mu+126\mu^2-42\mu^4).
\ee
From the above definition of $\gamma_j$ and (\ref{e46a}) we have
\bee\label{e46b}
\gamma_j=\alpha-j+1-\frac{1}{2\sigma}.
\ee

\vspace{0.4cm}

\noindent{\bf 4.2.\ The expansion ${}_0\Psi_1(x)$ as $x\to+\infty$ when $\sigma\neq\fs$}
\vspace{0.3cm}

\noindent
We now return to consideration of the expansion of the generalised Bessel function ${}_0\Psi_1(x)$ in (\ref{e35}) as $x\to+\infty$. The algebraic component of this expansion is, from (\ref{e35}) and (\ref{e37}) with $\delta=1-b$,
\begin{eqnarray}
{\hat H}_{0,1}(x)&:=&\frac{1}{2\pi}\{e^{\pi i\vartheta} H_{1,0}(xe^{\pi i\sigma}.\,e^{-\pi i})+e^{-\pi i\vartheta} H_{1,0}(xe^{-\pi i\sigma}.\,e^{\pi i})\}\nonumber\\
&=&\frac{1}{2\pi i\sigma}\sum_{k=0}^\infty \frac{\g(\frac{k+1-b}{\sigma})}{k!}\,\{(xe^{-\pi i})^{-(k+1-b)/\sigma}-(xe^{\pi i})^{-(k+1-b)/\sigma}\}\nonumber\\
&=&\frac{1}{\sigma} \sum_{k=0}^\infty \frac{x^{-(k+1-b)/\sigma}}{k! \g(1-\frac{k+1-b}{\sigma})}~.\label{e41}
\end{eqnarray}
In Section 4.1 we assumed that $\g(\sigma n+1-b)$ is regular for $n=0, 1, 2, \ldots\ $. If this assumption is false, the combination appearing in (\ref{e41}) is still valid; see \cite[\S 4]{W39}.  In such cases,  the expansion (\ref{e41}) becomes a finite sum when, in addition, $1/\sigma$ is an integer.

We define
\bee\label{e42d}
X_\pm=Xe^{\pm\pi i\sigma/\kappa},\qquad X=\kappa(hx)^{1/\kappa}.
\ee
Then the exponential component is 
\begin{eqnarray}
{\hat E}_{0,1}(x)\!\!\!&:=&\!\!\!\frac{1}{2\pi}\{e^{\pi i\vartheta}E_{0,1}(xe^{\pi i\sigma})+e^{-\pi i\vartheta} E_{0,1}(xe^{-\pi i\sigma})\}\nonumber\\
\!\!\!&=&\!\!\!\frac{1}{2\pi}\bl\{e^{\pi i\vartheta} X_+^\vartheta e^{X_+}\sum_{j=0}^\infty A_j(\sigma)X_+^{-j}+e^{-\pi i\vartheta} X_-^\vartheta e^{X_-}\sum_{j=0}^\infty A_j(\sigma)X_-^{-j}\br\}\nonumber\\
\!\!\!&=&\!\!\!\frac{X^\vartheta}{\pi}\,e^{X \cos \pi\sigma/\kappa} \sum_{j=0}^\infty (-)^jA_j(\sigma)X^{-j} \cos\,\bl[X \sin\frac{\pi\sigma}{\kappa}+\frac{\pi}{\kappa}(\vartheta\!-\!j)\br]\label{e42}
\end{eqnarray}
provided $0<\sigma<\fs$, where the first few coefficients $A_j(\sigma)$ can be obtained from (\ref{ecj}) with $\delta=1-b$ and the higher coefficients in specific cases by means of the algorithm in Appendix A.

Then we obtain the following theorem.
\begin{theorem}$\!\!\!.$\ When $a=-\sigma$, with $0<\sigma<1$, we have the expansions of the generalised Bessel function
\bee\label{e38}
{}_0\Psi_1(x)\sim\left\{\begin{array}{ll} {\hat E}_{0,1}(x)+{\hat H}_{0,1}(x) & (0<\sigma<\fs)\\
\\
{\hat H}_{0,1}(x) & (\fs<\sigma<1)\end{array}\right.
\ee
as $x\to+\infty$, where ${\hat H}_{0,1}(x)$ and ${\hat E}_{0,1}(x)$ are defined in (\ref{e41}) and (\ref{e42}). 
The coefficients $A_j(\sigma)$ can be obtained from (\ref{ecj}) with $\delta=1-b$.
\end{theorem}
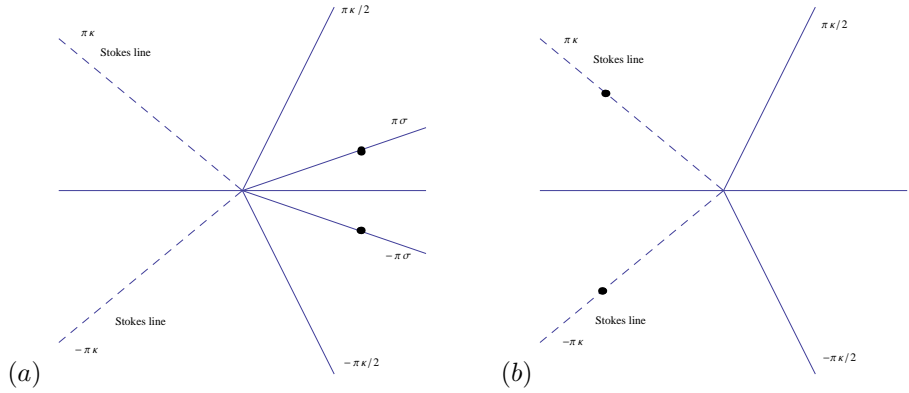
\begin{figure}[t]
\centering
\begin{picture}(300,200)(0,0)
\put(0,100){\line(1,0){120}}
\put(180,100){\line(1,0){120}}
\put(60,100){\line(-1,1){60}}
\put(60,100){\line(-1,-1){60}}
\put(60,100){\line(3,1){70}}
\put(60,100){\line(3,-1){70}}
\put(240,100){\line(-1,1){60}}
\put(240,100){\line(-1,-1){60}}
\multiput(60,100)(14,14){5}{\line(1,1){10}}
\multiput(60,100)(14,-14){5}{\line(1,-1){10}}
\multiput(240,100)(14,14){5}{\line(1,1){10}}
\multiput(240,100)(14,-14){5}{\line(1,-1){10}}
%First figure captions
\put(-10,166){\footnotesize{$\pi\kappa$}}
\put(-10,175){\footnotesize{Stokes line}}
\put(-10,30){\footnotesize{$-\pi\kappa$}}
\put(-10,20){\footnotesize{Stokes line}}
\put(90,166){\footnotesize{$\pi\kappa/2$}}
\put(85,26){\footnotesize{$-\pi\kappa/2$}}
\put(100,122){\footnotesize{$\pi\sigma$}}
\put(100,73){\footnotesize{$-\pi\sigma$}}
\put(-5,0){$(a)$}
%Second figure captions
\put(170,166){\footnotesize{$\pi\kappa$}}
\put(170,175){\footnotesize{Stokes line}}
\put(170,30){\footnotesize{$-\pi\kappa$}}
\put(170,20){\footnotesize{Stokes line}}
\put(270,166){\footnotesize{$\pi\kappa/2$}}
\put(265,26){\footnotesize{$-\pi\kappa/2$}}
\put(175,0){$(b)$}
\put(108,116){\circle*{4}}
\put(108,84){\circle*{4}}
\put(200,140){\circle*{4}}
\put(200,60){\circle*{4}}
\end{picture}
\caption{\small{The Stokes and anti-Stokes lines when $0<\kappa<1$ and the location of the arguments $P_\pm$ (indicated by heavy dots):
(a) $P_\pm=xe^{\pm\pi i\sigma}$ and (b) $P_\pm=xe^{\pm\pi i\kappa}$. }}
\end{figure}

Let us denote the points $xe^{\pm\pi i\sigma}$ in the $z$-plane that appear in the arguments of the associated function ${}_1\Psi_0(z)$ in (\ref{e35}) by $P_\pm$. Then when
$0<\sigma<\f{1}{3}$, $P_\pm$ lie in the exponentially large sector $|\arg\,z|<\fs\pi\kappa$ in Fig.~2(a) and consequently
${\hat E}_{0,1}(x)$ is exponentially large as $x\to+\infty$. When $\sigma=\f{1}{3}$, $P_\pm$ lie on the anti-Stokes lines $\arg\,z=\pm\fs\pi\kappa$; on these rays $\cos \pi\sigma/\kappa=0$ and ${\hat E}_{0,1}(x)$ is oscillatory with an algebraically controlled amplitude. When $\f{1}{3}<\sigma<\fs$, $P_\pm$ lie in the exponentially small and algebraic sectors so that ${\hat E}_{0,1}(x)$ is exponentially small. When $\sigma=\fs$, $P_\pm$ lie on the Stokes lines $\arg\,z=\pm\pi\kappa$, where the subdominant exponential expansions are in the process of switching off. Finally, when $\fs<\sigma<1$, $P_\pm$ are situated in the algebraic sectors, where the expansions are purely algebraic.
\vspace{0.4cm}

\noindent{\bf 4.3.\ The expansion ${}_0\Psi_1(x)$ as $x\to+\infty$ when $\sigma=\fs$}
\vspace{0.3cm}

\noindent
The expansion of ${}_0\Psi_1(x)$ when $\sigma=\fs$ requires a separate treatment.
In this case, the functions ${}_1\Psi_0(xe^{\pm\pi i\sigma})$ defined in (\ref{e36}) have arguments situated on the Stokes lines $\arg\,z=\pm\pi\kappa$, since $\kappa=1-\sigma=\fs$. 
Then, we obtain from (\ref{e46}) with $\delta=1-b$ the expansions
\[{}_1\Psi_0(xe^{\pm\fr\pi i})=\pm 2ie^{\pm\pi i\vartheta} \sum_{k=0}^{m_0-1} \frac{\g(2k+2-2b)}{k!} x^{-2(k+1-b)}\]
\bee\label{e400}
+(Xe^{\pm\pi i})^\vartheta e^{-X} \bl\{\sum_{j=0}^{M-1}\bl(\frac{1}{2}A_j(\fs)\pm \frac{i B_j(\fs)}{\sqrt{2\pi X}}\br)(-X)^{-j}+O(X^{-M})\br\}
\ee
as $x\to+\infty$, where $X=x^2/4$ and $\vartheta=\fs-b$. 

In Appendix B, we give a different derivation of (\ref{e46}) when $\sigma=\fs$ (where $\kappa=\fs$, $\mu=1$) by making use of the relation of ${}_1\Psi_0(z)$ to the confluent hypergeometric function. In this case the coefficients $A_j(\fs)$ appearing in (\ref{e400}) are given in closed form by
\bee\label{e42f}
A_j(\fs)=2\sqrt{\pi}\,\frac{(b)_j(b-\fs)_j}{j!};
\ee
see (\ref{b3}). The coefficients $B_j(\fs)$ are given by (\ref{e47}) with the  corresponding coefficients $G_{2k,j}(\sigma)$ (for $0\leq k\leq 4$) when $\sigma=\fs$ ($\mu=1$) in (\ref{e42e}) becoming\footnote{There was a misprint in the first term in ${\hat G}_{6,j}$ in \cite{PCHF}, which appeared as $-3226$ instead of $-3626$. This was pointed out by T. Pudlik \cite{TP}. The correct value was used in the numerical calculations described in \cite{PCHF}.}
\begin{eqnarray}
{\hat G}_{0,j}\!\!&=&\!\!\f{2}{3}-\gamma_j,\qquad {\hat G}_{2,j}=\f{1}{15}(46-225\gamma_j+270\gamma_j^2-90\gamma_j^3),\nonumber \\
{\hat G}_{4,j}\!\!&=&\!\!\f{1}{70}(230-3969\gamma_j+11340\gamma_j^2-11760\gamma_j^3+5040\gamma_j^4
-756\gamma_j^5),\nonumber\\
{\hat G}_{6,j}\!\!&=&\!\!\f{1}{350}(-3626-17781\gamma_j+183330\gamma_j^2-397530\gamma_j^3+370440\gamma_j^4
-170100\gamma_j^5\nonumber\\
&&\hspace{7cm}+37800\gamma_j^6-3240\gamma_j^7),\nonumber\\
{\hat G}_{8,j}\!\!&=&\!\!\f{1}{231000}(-4032746+43924815\gamma_j+88280280\gamma_j^2-743046480\gamma_j^3\nonumber\\
&&+1353607200\gamma_j^4-1160830440\gamma_j^5+541870560\gamma_j^6
-141134400\gamma_j^7\nonumber\\
&&\hspace{6cm}+19245600\gamma_j^8-1069200\gamma_j^9),\label{e42c}
\end{eqnarray}
where ${\hat G}_{2k,j}=6^{2k} G_{2k,j}(\fs)$ and, from (\ref{e46b}), $\gamma_j=\alpha-j$. Higher coefficients can be obtained from (\ref{e49}) and (\ref{e49a}).

Insertion of the expansion (\ref{e400}) into (\ref{e35}) with $\delta=1-b$, followed by some routine algebra, then yields the following
\begin{theorem}$\!\!\!.$\ Let $\vartheta=\fs-b$, $X=x^2/4$ and $M$ denote a positive integer. Then, when $\sigma=\fs$ we have the expansion for the generalised Bessel function
\[{}_0\Psi_1(x)={\hat H}_{0,1}^{opt}(x)+\frac{X^\vartheta\,e^{-X}}{2\pi}\bl\{ \sum_{j=0}^{M-1}\bl( \cos 2\pi\vartheta\,A_j(\fs)-\frac{2\sin 2\pi\vartheta\,B_j(\fs)}{\sqrt{2\pi X}}\br)(-X)^{-j}\]
\bee\label{e45}
\hspace{7cm}+O(X^{-M})\br\}
\ee
as $x\to+\infty$, where ${\hat H}_{0,1}^{opt}(x)$ is the optimally truncated algebraic expansion obtained from (\ref{e41})
\bee\label{e45a}
{\hat H}_{0,1}^{opt}(x)=2\sum_{k=0}^{m_o-1} \frac{x^{-2(k+1-b)}}{k!\,\g(2b-1-2k))},\qquad m_o=X+2b-\f{3}{2}+\alpha
\ee
with $|\alpha|$ bounded. The coefficients $A_j(\fs)$ and $B_j(\fs)$ are given in (\ref{e42f}) and (\ref{e47}), where
$G_{2k,j}(\fs)$ are those appearing in (\ref{e42c}) with $\gamma_j=\alpha-j$, $|\alpha|$ bounded.     
\end{theorem} 
We remark that when $2b=1, 0, -1, -2, \ldots$ the algebraic expansion in (\ref{e45a}) vanishes and that when $2b=2, 3 \ldots$ it reduces to a finite sum. In both these situations the parameter $2\vartheta$ is an integer and so the second sum in the exponentially small contribution in (\ref{e45}) is not present.

\vspace{0.4cm}

\noindent{\bf 4.4.\ The expansion of ${}_0\Psi_1(-x)$ as $x\to+\infty$}
\vspace{0.3cm}

\noindent 
When $z=-x$, $x>0$ we find from (\ref{e35}) and (\ref{e36}), upon replacing $z$ by $e^{\mp\pi i}x$ and using the result ${}_1\Psi_0(ze^{2\pi i})={}_1\Psi_0(z)$, that
\bee\label{e39}
{}_0\Psi_1(-x)=\frac{1}{2\pi}\bl\{e^{\pi i\vartheta} {}_1\Psi_0(xe^{-\pi i\kappa})+e^{-\pi i\vartheta}{}_1\Psi_0(xe^{\pi i\kappa})\br\}.
\ee
It is now apparent that the arguments of the associated ${}_1\Psi_0$ functions are situated on the Stokes lines for all values of $\sigma$ in the range $0<\sigma<1$, where the exponential expansions will be in the process of switching off; see Fig.~2(b). 

The algebraic component of the right-hand side of (\ref{e39}) is, from (\ref{e46}),
\[2\cos \pi(\delta-\vartheta) \sum_{k=0}^{m_o-1} \frac{\g(\frac{k+\delta}{\sigma})}{k!}\,x^{-(k+\delta)/\sigma} \equiv 0,\]
upon recalling that $\delta=\vartheta+\fs$.
The exponentially small contributions involving the coefficients $B_j(\sigma)$ in (\ref{e46}) are seen to cancel 
in the combination in (\ref{e39}), thereby yielding the final result:
\begin{theorem}$\!\!\!.$\ Let $\vartheta=\fs-b$ and $M$ denote a positive integer. Then, when $0<\sigma<1$ we have the expansion of the generalised Bessel function 
\bee\label{e410}
{}_0\Psi_1(-x)=\frac{X^\vartheta\,e^{-X}}{2\pi} \bl\{\sum_{j=0}^{M-1} (-)^jA_j(\sigma)X^{-j}+O(X^{-M})\br\} \qquad (0<\sigma<1)
\ee
as $x\to+\infty$, where $X$ is defined in (\ref{e42d}). The first few coefficients $A_j(\sigma)$ can be obtained from (\ref{ecj}) with $A_0(\sigma)$ given in (\ref{e37a}).
\end{theorem}
\vspace{0.6cm}

%\begin{center}
{\bf 5. \ The expansion of ${}_0\Psi_1(z)$ when $-1<a<0$ for $|z|\to\infty$}
%\end{center}
\setcounter{section}{5}
\setcounter{equation}{0}
\renewcommand{\theequation}{\arabic{section}.\arabic{equation}}
From (\ref{e35}), the generalised Bessel function when $a=-\sigma$, $0<\sigma<1$, is given by
\bee\label{e51}
{}_0\Psi_1(z)=\frac{1}{2\pi}\bl\{e^{\pi i\vartheta} {}_1\Psi_0(ze^{\pi i\sigma})+e^{-\pi i\vartheta} {}_1\Psi_0(ze^{-\pi i\sigma})\br\},
\ee
where the associated function ${}_1\Psi_0(z)$ is defined in (\ref{e36}) with $\delta=1-b=\fs+\vartheta$. The division of the $z$-plane where ${}_1\Psi_0(z)$ possesses exponentially large and algebraic behaviour as $|z|\to\infty$, and the location of the Stokes lines on the rays $\arg\,z=\pm\pi\kappa$, is discussed in Section 4.2.

In addition to the above-mentioned Stokes lines, the positive real axis $\arg\,z=0$ is also a Stokes line for the algebraic expansion. On the positive real axis, where the algebraic expansion in (\ref{e37}) is maximally subdominant, the leading coefficient of $z^{-\delta/\sigma}$ multiplying the algebraic expansion changes by the factor $e^{2\pi i\delta/\sigma}$ as $\arg\,z$ changes from small negative to small positive values.  The result of this change is that in the combination (\ref{e51}) when $\theta\pm\pi\sigma$ are of different signs, the algebraic expansions $H_{1,0}(ze^{\pm\pi i\sigma})$ combine to yield the expansion ${\hat H}_{0,1}(z)$ (compare (\ref{e41})), since
\[\frac{1}{2\pi}\bl\{e^{\pi i\vartheta} H_{1,0}(ze^{\pi i\sigma} .\, e^{-\pi i})+e^{-\pi i\vartheta} H_{1,0}(ze^{-\pi i\sigma} . \,e^{\pi i})\br\}\]
\bee\label{e52}
=\frac{1}{\pi\sigma}\sum_{k=0}^\infty \frac{\g(\frac{k+\delta}{\sigma})\sin\,\pi(\frac{k+\delta}{\sigma})}{k!} \,z^{-(k+\delta)/\sigma}= \frac{1}{\sigma} \sum_{k=0}^\infty \frac{z^{-(k+\delta)/\sigma}}{k! \g(1-\frac{k+\delta}{\sigma})}\,.
\ee
However, when $\theta\pm\pi\sigma$ are of the same sign, the algebraic expansions in the combination (\ref{e21}) cancel, since  we have for $\theta\pm\pi\sigma \gl 0$
\[e^{\pi i\vartheta} H_{1,0}(ze^{\pi i\sigma} .\, e^{\mp\pi i})+e^{-\pi i\vartheta} H_{1,0}(ze^{-\pi i\sigma} . \,e^{\mp\pi i})\]
\[=2\cos \pi(\delta-\vartheta)\sum_{k=0}^\infty \frac{\g(\frac{k+\delta}{\sigma})}{k!} \,(ze^{\mp\pi i})^{-(k+\delta)/\sigma}\equiv 0.\]
\begin{figure}[t]
\centering
\begin{picture}(400,200)(0,0)
\put(0,100){\line(1,0){100}}
\put(150,100){\line(1,0){100}}
\put(300,100){\line(1,0){100}}
\put(50,100){\line(-1,1){50}}
\put(50,100){\line(-1,-1){50}}
\put(50,100){\line(3,1){60}}
\put(50,100){\line(3,-1){60}}
\put(200,100){\line(-1,1){50}}
\put(200,100){\line(-1,-1){50}}
\put(200,100){\line(1,2){30}}
\put(200,100){\line(1,-2){30}}
\put(350,100){\line(-1,1){50}}
\put(350,100){\line(-1,-1){50}}
\put(350,100){\line(1,6){10}}
\put(350,100){\line(1,-6){10}}
\multiput(50,100)(14,14){4}{\line(1,1){10}}
\multiput(50,100)(14,-14){4}{\line(1,-1){10}}
\multiput(200,100)(14,14){4}{\line(1,1){10}}
\multiput(200,100)(14,-14){4}{\line(1,-1){10}}
\multiput(350,100)(14,14){4}{\line(1,1){10}}
\multiput(350,100)(14,-14){4}{\line(1,-1){10}}
%First figure labels
\put(-10,166){\footnotesize{$\pi\kappa$}}
\put(-10,175){\footnotesize{Stokes line}}
\put(-10,30){\footnotesize{$-\pi\kappa$}}
\put(-10,20){\footnotesize{Stokes line}}
\put(90,166){\footnotesize{$\pi\kappa/2$}}
\put(85,26){\footnotesize{$-\pi\kappa/2$}}
\put(100,122){\footnotesize{$\pi\sigma$}}
\put(100,73){\footnotesize{$-\pi\sigma$}}
\put(-5,0){$(a)$}
%Second figure labels
\put(150,166){\footnotesize{$\pi\kappa$}}
\put(150,175){\footnotesize{Stokes line}}
\put(150,30){\footnotesize{$-\pi\kappa$}}
\put(150,20){\footnotesize{Stokes line}}
\put(250,166){\footnotesize{$\pi\kappa/2$}}
\put(245,26){\footnotesize{$-\pi\kappa/2$}}
\put(220,166){\footnotesize{$\pi\sigma$}}
\put(215,26){\footnotesize{$-\pi\sigma$}}
\put(150,0){$(b)$}
%Third figure labels
\put(355,166){\footnotesize{$\pi\kappa$}}
\put(355,175){\footnotesize{Stokes line}}
\put(355,30){\footnotesize{$-\pi\kappa$}}
\put(355,20){\footnotesize{Stokes line}}
\put(390,130){\footnotesize{$\pi\kappa/2$}}
\put(385,70){\footnotesize{$-\pi\kappa/2$}}
\put(315,140){\footnotesize{$\pi\sigma$}}
\put(310,55){\footnotesize{$-\pi\sigma$}}
\put(300,0){$(c)$}

\put(98,116){\circle*{4}}
\put(98,84){\circle*{4}}
\put(225,150){\circle*{4}}
\put(225,50){\circle*{4}}
\put(320,130){\circle*{4}}
\put(320,70){\circle*{4}}
\end{picture}
\caption{\small{The Stokes and anti-Stokes lines and the location of the arguments $xe^{\pm\pi i\sigma}$ (indicated by heavy dots) when
(a) $0<\sigma<\f{1}{3}$,  (b) $\f{1}{3}<\sigma<\fs$ and (c) $\fs<\sigma<1$. }}
\end{figure}
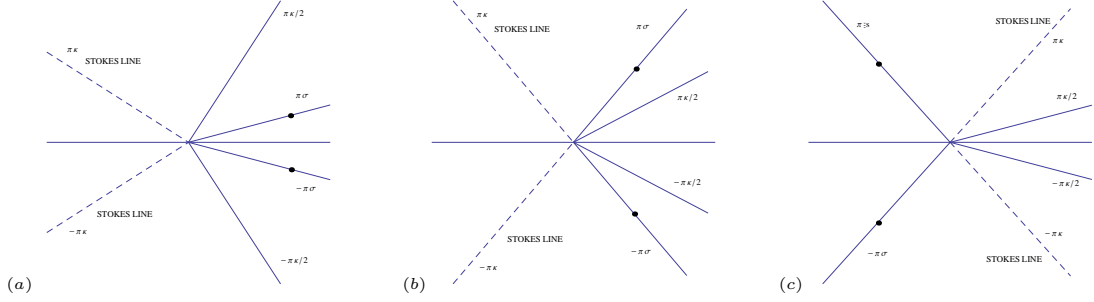

The expansion of the generalised Bessel function ${}_0\Psi_1(z)$ with $0<\sigma<1$ can then be constructed from (\ref{e51}) by application of the expansions of ${}_1\Psi_0(z)$ in Section 4.1. We consider separately the ranges (i) $0<\sigma<\f{1}{3}$, (ii) $\sigma=\f{1}{3}$, (iii) $\f{1}{3}<\sigma\leq\fs$ and (iv) $\fs<\sigma<1$. In most cases we omit the values $\theta=0$ and $\theta=\pm\pi$ as these are covered by Theorems  4 and 6. For convenience in presentation we define the exponential expansions
\begin{eqnarray}
E_\pm(z)&:=&\frac{e^{\pm\pi i\vartheta}}{2\pi}\,E_{1,0}(ze^{\pm\pi i\sigma})\nonumber\\
&=&\frac{(e^{\pm\pi i} Z_\pm)^\vartheta}{2\pi}\,
e^{Z_\pm} \sum_{j=0}^\infty A_j(\sigma) Z_\pm^{-j},\qquad Z_\pm=X e^{i(\theta\pm\pi\sigma)/\kappa},\label{e26}
\end{eqnarray}
where $X$ is defined in (\ref{e42d}).
The algebraic expansion ${\hat H}_{0,1}(z)$ is defined in (\ref{e52}).
\bigskip

\noindent (i)\ Case $0<\sigma<\f{1}{3}$:\ \ \ Fig.~3(a) shows the positions of the points $ze^{\pm\pi i\sigma}$ when $\theta=0$ and $0<\sigma<\f{1}{3}$ in relation to the rays $\pm\fs\pi\kappa$ and $\pm\pi\kappa$. Then taking into account the Stokes lines on the rays $\pm\pi\kappa$ and the positive real axis, we have from (\ref{e51}) and (\ref{e4exp})
\bee\label{e27e}
{}_0\Psi_1(z)\sim\left\{\begin{array}{ll} E_+(z)+E_-(z)+{\hat H}_{0,1}(z) & (|\theta|<\pi\sigma)\\
\\
E_+(z)+E_-(z) & (\pi\sigma<|\theta|<\pi(1-2\sigma))\\
\\
E_\mp(z) & (\pi(1-2\sigma)<|\theta|<\pi)\\
\end{array}\right.
\ee
as $|z|\to\infty$. In the third expansion the upper or lower sign is chosen according as $\theta>0$ or $\theta<0$, respectively. If the subdominant algebraic expansion present in the sector $|\theta|<\pi\sigma$ is neglected, then the expansion $E_+(z)+E_-(z)$ holds in the wider sector $|\theta|<\pi(1-2\sigma)$. 
\bigskip

\noindent (ii)\ Case $\sigma=\f{1}{3}$:\ \ \ When $\sigma=\f{1}{3}$ and $\theta=0$ the points $ze^{\pm\pi i\sigma}$ are situated on the anti-Stokes rays $\pm\f{1}{3}\pi$. Then taking into account the Stokes lines on the rays $\pm\f{2}{3}\pi$ and the positive real axis, we have from (\ref{e51}) and (\ref{e4exp})
\bee\label{e27ae}
{}_0\Psi_1(z)\sim\left\{\begin{array}{ll} E_+(z)+E_-(z)+{\hat H}_{0,1}(z) & (|\theta|<\f{1}{3}\pi)\\
\\
E_\mp(z) & (\f{1}{3}\pi\leq|\theta|<\pi)\\
\end{array}\right.
\ee
as $|z|\to\infty$. 
When $\theta=\pm\f{1}{3}\pi$, the points $ze^{\pm\pi i\sigma}$ are situated either on one of the rays $\pm\f{2}{3}\pi$, where the exponential expansions $E_\pm(z)$ are in the process of switching off, or the positive real axis, where the subdominant algebraic expansion is also in the process of switching off. Thus we have a double Stokes phenomenon when $\theta=\pm\f{1}{3}\pi$. We do not consider this case in detail, but since $E_\mp(z)$ will be maximally dominant when $\theta=\pm\f{1}{3}\pi$, we can say that, neglecting the subdominant algebraic expansion and the sub-subdominant exponential expansions $E_\pm(z)$, the expansion of ${}_0\Psi_1(z)$ for large $|z|$ when $\theta=\pm\f{1}{3}\pi$ is given by $E_\mp(z)$.
\bigskip

\noindent (iii)\ Case $\f{1}{3}<\sigma\leq\fs$:\ \ \ Fig.~3(b) shows the situation when $\f{1}{3}<\sigma<\fs$; when $\sigma=\fs$ the points $ze^{\pm\pi i\sigma}$ when $\theta=0$ lie on the rays $\pm\pi\kappa$ (where $\kappa=\fs)$.
Then taking into account the Stokes lines on the rays $\pm\pi\kappa$ and the positive real axis, we have from (\ref{e51}) and (\ref{e4exp})
\bee\label{e28e}
{}_0\Psi_1(z)\sim\left\{\begin{array}{ll} E_+(z)+E_-(z)+{\hat H}_{0,1}(z) & (|\theta|<\pi(1-2\sigma))\\
\\
E_\mp(z)+{\hat H}_{0,1}(z) & (\pi(1-2\sigma)<|\theta|<\pi\sigma)\\
\\
E_\mp(z) & (\pi\sigma<|\theta|<\pi)\\
\end{array}\right.
\ee
as $|z|\to\infty$. When $\sigma=\fs$ the first expansion in (\ref{e28e}) is inapplicable.
\bigskip

\noindent (iv)\ Case $\fs<\sigma<1$:\ \ \ Fig.~3(c) shows the situation when $\fs<\sigma<1$. The points $ze^{\pm\pi i\sigma}$ when $\theta=0$ lie in the algebraic sector. 
Then taking into account the Stokes lines on the rays $\pm\pi\kappa$ and the positive real axis, we have from (\ref{e51}) and (\ref{e4exp})

\bee\label{e29e}
{}_0\Psi_1(z)\sim\left\{\begin{array}{ll} {\hat H}_{0,1}(z) & (|\theta|<\pi(2\sigma-1))\\
\\
E_\mp(z)+{\hat H}_{0,1}(z) & (\pi(2\sigma-1)<|\theta|<\pi\sigma)\\
\\
E_\mp(z) & (\pi\sigma<|\theta|<\pi)\\
\end{array}\right.
\ee
as $|z|\to\infty$. We note that the algebraic expansion in the middle expansion in (\ref{e29e}) is dominant in the sectors $\pi(2\sigma-1)<|\theta|<\pi(3\sigma-1)/2$. Hence if the subdominant exponential expansion is neglected, then the expansion ${\hat H}_{0,1}(z)$ holds in the wider sector $|\theta|<\pi(3\sigma-1)/2$.

%\begin{center}
{\bf 6. \  Numerical results}
%\end{center}
\setcounter{section}{6}
\setcounter{equation}{0}
\renewcommand{\theequation}{\arabic{section}.\arabic{equation}}
We present some numerical results to verify the statements in Theorems 5 and 6. For the generalised Bessel function\footnote{For convenience in this section we denote the dependence of ${}_0\Psi_1(x)$ on $\sigma$ by writing ${}_0\Psi_1(\sigma;x)$.}
with $\sigma=\fs$
\[{}_0\Psi_1(x)\equiv{}_0\Psi_1(\fs;x)=\sum_{n=0}^\infty \frac{x^n}{n!\,\g(-\fs n+b)},\]
we first observe that when $b=-m+\fs$ or $b=-m$, $m=0, 1, 2, \ldots $ we have the exact representations
\bee\label{e61}
{}_0\Psi_1(\fs;x)=\pm \frac{X^\vartheta e^{-X}}{\sqrt{\pi}} \sum_{j=0}^m (-m\pm\fs)_j\bl(\!\!\!\begin{array}{c}m\\j\end{array}\!\!\!\br) X^{-j},\qquad b=\left\{\begin{array}{l}\!\!-m+\fs\\ \!\!-m,\end{array}\right.
\ee
where $X=x^2/4$, $\vartheta=\fs-b$ and the upper (resp. lower) signs correspond to $b=-m+\fs$ (resp. $b=-m$). This result can be established by use of the confluent hypergeometric function representation in (\ref{b1}), together with Kummer's transformation to extract the factor $e^{-X}$, followed by some routine algebra. With these values of $b$, the algebraic expansion in (\ref{e45a}) vanishes identically and $\sin 2\pi\vartheta=0$. The coefficients $A_j(\fs)$ in (\ref{e42f}) vanish for $j\geq m+1$, whence it follows that the expansions in Theorems 5 and 6 (when $\sigma=\fs$) agree with (\ref{e61}).

We now consider other values of the parameter $b$ not satisfying those in (\ref{e51}).
From Theorem 5, we have when $\sigma=\fs$ and $X=x^2/4$
\bee\label{e62}
{}_0\Psi_1(\fs;x)-H_{0,1}^{opt}(x)\sim R_M^+(x)\qquad (x\to+\infty),
\ee
where
\[R_M^+(x):=\frac{X^\vartheta e^{-X}}{2\pi}\sum_{j=0}^{M}\bl\{ \cos 2\pi\vartheta\,A_j(\fs)-\frac{2\sin 2\pi\vartheta\,B_j(\fs)}{\sqrt{2\pi X}}\br\}(-X)^{-j}.\]
The coefficients $B_j(\fs)$ are computed from (\ref{e47}).
In Table 2 we compare  the values of the left- and right-hand sides of (\ref{e62}) for different values of $b$ and truncation index $M$ when $x=10$. 

The first value, $b=1$, corresponds to a terminating algebraic expansion (which terminates at the leading term equal to 2). In this case, the quantity $2\vartheta$ is an integer so that the component of $R_M^+(x)$ involving the $B_j(\fs)$ is absent. We have optimally truncated the exponential expansion at $M=24$. The second value $b=\f{1}{4}$ (corresponding to $m_o=24$, $\alpha=0$) yields an $R_M^+(x)$ in which the first component involving the $A_j(\fs)$ is absent; here we have taken $M=6$. The remaining values in the table correspond to an optimally truncated algebraic expansion and $R_M^+(x)$ involving both component expansions (with $M=6$). In the case $b=\f{1}{3}$, we have from (\ref{e45a}) $m_o=24$ and $\alpha=-\f{1}{6}$. In \cite[\S 5]{W39}, Wright gave the expansion (\ref{e62}) for ${}_0\Psi_1(x)$ when $\sigma=\fs$, but with the exponentially small contribution given by (in our notation)
\[\frac{X^\vartheta e^{-X}}{2\pi} \sum_{j=0}^\infty A_j(\fs) (-X)^{-j}\]
for $x\to+\infty$.
As our results demonstrate this is, in general, not correct as it does not take into account the Stokes phenomenon.
\begin{table}[th]
\caption{\footnotesize{Values of ${}_0\Psi_1(\fs;x)-H_{0,1}^{opt}(x)$ and $R_M^+(x)$ for different values of $b$ when $x=10$. The truncation index $M$ is indicated in the text. }}
\begin{center}
\begin{tabular}{|c|l|l|}
\hline
&& \\[-0.4cm]
\mcol{1}{|c|}{$b$} & \mcol{1}{c|}{${}_0\Psi_1(\fs;x)-H_{0,1}^{opt}(x)$} & \mcol{1}{c|}{$R_M^+(x)$}\\
[.05cm]\hline
%&& \\[-0.45cm]
%\fs & $+7.8354332655\times 10^{-12}$ & $+7.8354332655\times 10^{-12}$ \\
&& \\[-0.35cm]
1   & $-1.5374597944\times 10^{-12}$ & $-1.5374597944\times 10^{-12}$  \\
&& \\[-0.4cm]
$\f{1}{4}$ & $-1.8851189300\times 10^{-12}$ & $-1.8851189236\times 10^{-12}$ \\
&& \\[-0.4cm]
$\f{1}{3}$ & $+5.1505426736\times 10^{-12}$ & $+5.1505426725\times 10^{-12}$ \\
&& \\[-0.4cm]
$\f{4}{5}$ & $-5.7125964076\times 10^{-13}$ & $-5.7125962370\times 10^{-13}$ \\
&& \\[-0.4cm]
$\f{6}{5}$ & $-3.1143753823\times 10^{-13}$ & $-3.1143783601\times 10^{-13}$ \\
[0.15cm]
\hline
\end{tabular}
\end{center}
\end{table}

In Table 3 we compare the values of ${}_0\Psi_1(\sigma;-x)$ and $R_M^-(x)$ for different values of $\sigma$ when $b=1$ and $x=10$, where from Theorem 6
\[{}_0\Psi_1(\sigma;-x)\sim R_M^-(x),\qquad R_M^-(x):=\frac{X^\vartheta e^{-X}}{2\pi}\sum_{j=0}^{M}A_j(\sigma) (-X)^{-j}.\]
The first few coefficients $A_j(\sigma)$ are obtained from (\ref{ecj}), with $A_0(\sigma)$ given in (\ref{e37a}), and higher coefficients from the algorithm in Appendix A. In all cases we have employed the truncation index $M=15$.
\begin{table}[th]
\caption{\footnotesize{Values of ${}_0\Psi_1(\sigma;-x)$ and $R_M^-(x)$ for different values of $\sigma$ and $x$ when $b=1$. The truncation index $M=15$. }}
\begin{center}
\begin{tabular}{|c|c|l|l|}
\hline
&& \\[-0.4cm]
\mcol{1}{|c|}{$\sigma$} & \mcol{1}{c|}{$x$} & \mcol{1}{c|}{${}_0\Psi_1(\sigma;-x)$} & \mcol{1}{c|}{$R_M^-(x)$}\\
[.05cm]\hline
%&& \\[-0.45cm]
%\fs & $+7.8354332655\times 10^{-12}$ & $+7.8354332655\times 10^{-12}$ \\
&&& \\[-0.35cm]
$\f{1}{4}$ & 15 & $+4.7317589195\times 10^{-9}$  & $+4.7317587800\times 10^{-9}$  \\
&&& \\[-0.4cm]
$\f{1}{3}$ & 12 & $+1.8807037460\times 10^{-8}$  & $+1.8807035571\times 10^{-8}$ \\
%&&& \\[-0.4cm]
%$\f{1}{3}$ & 15 & $+2.7740903726\times 10^{-11}$ & $+2.7740903971\times 10^{-11}$ \\
&&& \\[-0.4cm]
$\f{1}{2}$ & 10 & $-1.5374597944\times 10^{-12}$ & $-1.5374597943\times 10^{-12}$ \\
&&& \\[-0.4cm]
$\f{2}{3}$ & 6  & $+1.0783972342\times 10^{-15}$ & $+1.0783972342\times 10^{-15}$ \\
&&& \\[-0.4cm]
$\f{3}{4}$ & 4  & $+1.6389907960\times 10^{-13}$ & $+1.6389907960\times 10^{-13}$ \\
[0.15cm]
\hline
\end{tabular}
\end{center}
\end{table}

In the final part of this section, 
we carry out numerical computations to verify the expansions in (\ref{e27e})--(\ref{e29e}). For each value of $\sigma$ chosen it is necessary to compute the coefficients $A_j(\sigma)$ that appear in the exponential expansions $E_\pm(z)$.
The first three coefficients can be obtained from (\ref{ecj}), with higher coefficients generated by means of the algorithm in Appendix A. 
Examples of the coefficients for several values of $\sigma$ are given in Table 4; these values are employed in the calculations presented in Tables 5 and 6.
\begin{table}[th]
\caption{\footnotesize{The normalised coefficients $c_j(\sigma)=A_j(\sigma)/A_0(\sigma)$ for $1\leq j\leq 10$ when $b=\f{5}{4}$. }}
\begin{center}
\begin{tabular}{|c|l|l|l|}
\hline
\mcol{1}{|c|}{$j$} & \mcol{1}{c|}{$\sigma=1/6$} & \mcol{1}{c|}{$\sigma=1/3$} & \mcol{1}{c|}{$\sigma=2/3$}\\
[.05cm]\hline
&&& \\[-0.4cm]
1 & $1.86805555556$ & $1.16319444444$ & $0.83159722222$\\
2 & $5.71703800154$ & $2.59491343557$ & $1.53740023389$\\
3 & $2.32181131692\times 10^1$ & $8.42530200402$ & $4.38966463732$\\
4 & $1.16570408563\times 10^2$ & $3.58179860428\times 10^1$ & $1.69388501423\times 10^1$\\
5 & $6.98089732047\times 10^2$ & $1.88123659351\times 10^2$ & $8.23410445252\times 10^1$\\
6 & $4.87231305227\times 10^3$ & $1.17617708621\times 10^3$ & $4.82738754544\times 10^2$\\
7 & $3.89191967771\times 10^4$ & $8.52942466133\times 10^3$ & $3.31345555254\times 10^3$\\
8 & $3.50286638479\times 10^5$ & $7.03803279143\times 10^4$ & $2.60596676873\times 10^4$\\
9 & $3.50538397688\times 10^6$ & $6.51101116490\times 10^5$ & $2.31033323525\times 10^5$\\
10& $3.85813836005\times 10^7$ & $6.67440397372\times 10^6$ & $2.27941435603\times 10^6$\\
[.15cm]\hline
\end{tabular}
\end{center}
\end{table}

To verify the various expansions, we consider only $\theta\geq 0$ and define the quantity
\[F_{\lambda,\mu,\nu}(z):=\bl|\frac{\la E_+(z)+\mu E_-(z)+\nu {\hat H}_{1,0}(z)}{{}_0\Psi_1(z)}-1\br|\]
representing the absolute relative error in the computation of ${}_0\Psi_1(z)$,
where\footnote{The parameter $\mu$ here is not to be confused with that appearing in (\ref{e49}).} $\la$, $\mu$, $\nu=0$ or 1. This enables us to either retain or switch off the expansions $E_\pm(z)$ and ${\hat H}_{1,0}(z)$ by making an appropriate choice for $(\la,\mu,\nu)$. In order to detect the algebraic expansion when it is subdominant the exponential expansions were optimally truncated. For the values of $x$ in Tables 5 and 6 this required at most about 20 coefficients $A_j(\sigma)$. Similarly, in Table 6, when $\sigma=\f{2}{3}$ and $|\theta|<\f{1}{3}\pi$, the algebraic expansion has been optimally truncated.

\begin{table}[th]
\caption{\footnotesize{Values of $F_{\lambda,\mu,\nu}(z)$ for different values of $\theta$ and $\sigma$ when $z=xe^{i\theta}$ and $b=\f{5}{4}$. }}
\begin{center}
\begin{tabular}{|l|c|l||l|c|l|}
\hline
\mcol{3}{|c||}{$\sigma=1/6,\ x=15$} & \mcol{3}{c|}{$\sigma=1/3,\ x=10$}\\
%&&&&&&\\[-0.35cm]
\mcol{1}{|c|}{$\theta/\pi$} & \mcol{1}{c|}{$(\la, \mu, \nu)$} & \mcol{1}{c||}{$F_{\lambda,\mu,\nu}(z)$} & \mcol{1}{c|}{$\theta/\pi$} & \mcol{1}{c|}{$(\la, \mu, \nu)$} & \mcol{1}{c|}{$F_{\lambda,\mu,\nu}(z)$}\\
[.05cm]\hline
&&&&&\\[-0.35cm]
0 &    (1,1,1) & $3.422\times 10^{-6}$ & 0 &    (1,1,1) & $1.043\times 10^{-8}$\\ 
  &    (1,1,0) & $3.519\times 10^{-3}$ & &      (1,1,0) & $1.009\times 10^{0}$\\
&&&&&\\[-0.4cm]
0.10 & (1,1,1) & $9.554\times 10^{-5}$ & 0.10 & (1,1,1) & $6.426\times 10^{-5}$\\
  &    (1,1,0) & $5.199\times 10^{-4}$ & &      (1,1,0) & $4.797\times 10^{-1}$\\
&&&&&\\[-0.4cm]
0.25 & (1,1,0) & $8.495\times 10^{-5}$ & 0.25 & (1,1,1) & $1.504\times 10^{-4}$\\
  &    (1,1,1) & $6.700\times 10^{-4}$& &       (1,1,0) & $1.588\times 10^{-3}$\\
&&&&&\\[-0.4cm]
0.50 & (1,1,0) & $3.400\times 10^{-5}$ & 0.40 & (0,1,0) & $1.702\times 10^{-4}$\\
  &    (1,1,1) & $1.294\times 10^{+1}$ & &      (0,1,1) & $1.070\times 10^{-3}$\\
&&&&&\\[-0.4cm]
0.75 & (0,1,0) & $2.534\times 10^{-5}$ & 0.60 & (0,1,0) & $6.872\times 10^{-5}$\\
  &    (1,1,0) & $4.135\times 10^{-3}$ & &      (0,1,1) & $3.328\times 10^{+1}$\\
&&&&&\\[-0.4cm]
0.80 & (0,1,0) & $2.462\times 10^{-5}$ & 0.80 & (0,1,0) & $5.029\times 10^{-5}$\\
  &    (1,1,0) & $1.081\times 10^{-1}$ & &      (0,1,1) & $1.974\times 10^{+5}$\\      
[0.15cm]
\hline
\end{tabular}
\end{center}
\end{table}
\begin{table}[th]
\caption{\footnotesize{Values of $F_{\lambda,\mu,\nu}(z)$ for different values of $\theta$ and $\sigma$ when $z=xe^{i\theta}$ and $b=\f{5}{4}$. }}
\begin{center}
\begin{tabular}{|l|c|l||l|c|l|}
\hline
\mcol{3}{|c||}{$\sigma=2/5,\ x=10$} & \mcol{3}{c|}{$\sigma=2/3,\ x=5$}\\
%&&&&&&\\[-0.35cm]
\mcol{1}{|c|}{$\theta/\pi$} & \mcol{1}{c|}{$(\la,\mu,\nu)$} & \mcol{1}{c||}{$F_{\lambda,\mu,\nu}(z)$} & \mcol{1}{c|}{$\theta/\pi$} & \mcol{1}{c|}{$(\la,\mu,\nu)$} & \mcol{1}{c|}{$F_{\lambda,\mu,\nu}(z)$}\\
[.05cm]\hline
&&&&&\\[-0.35cm]
0 &    (1,1,1) & $1.171\times 10^{-10}$ & 0 &    (0,0,1) & $3.301\times 10^{-12}$\\ 
  &    (1,1,0) & $1.000\times 10^{0}$ & &        (0,1,1) & $1.741\times 10^{+6}$\\
&&&&&\\[-0.4cm]
0.10 & (1,1,1) & $1.857\times 10^{-8}$ & 0.10 & (0,0,1) & $4.135\times 10^{-12}$\\
  &    (1,1,0) & $0.995\times 10^{0}$ & &      (0,1,1) & $8.218\times 10^{+2}$\\
&&&&&\\[-0.4cm]
0.30 & (0,1,1) & $5.526\times 10^{-6}$ & 0.40 & (0,1,1) & $1.730\times 10^{-11}$\\
  &    (0,1,0) & $2.451\times 10^{-4}$& &       (0,0,1) & $4.464\times 10^{-9}$\\
&&&&&\\[-0.4cm]
0.50 & (0,1,0) & $5.527\times 10^{-6}$ & 0.50 & (0,1,1) & $5.886\times 10^{-10}$\\
  &    (0,1,1) & $2.451\times 10^{-4}$ & &      (0,0,1) & $1.499\times 10^{-2}$\\
&&&&&\\[-0.4cm]
0.60 & (0,1,0) & $3.313\times 10^{-6}$ & 0.80 & (0,1,0) & $4.659\times 10^{-8}$\\
  &    (1,1,0) & $1.000\times 10^{0}$ & &      (0,1,1) & $2.187\times 10^{-1}$\\
%&&&&&\\[-0.4cm]
%0.80 & 0 & 0 & $2.462\times 10^{-5}$ & 0.80 & 0 & 0 & $5.029\times 10^{-5}$\\
%  &    1 & 0 & $1.081\times 10^{-1}$ & &      0 & 1 & $1.974\times 10^{+5}$\\      
[0.15cm]
\hline
\end{tabular}
\end{center}
\end{table}

The first entry in Table 5 shows $\sigma=\f{1}{6}$ when the expansions in (\ref{e27e}) apply. It is seen that for $\theta=0$ and $0.1\pi$ (inside the sector $|\theta|<\f{1}{6}\pi$) the absolute relative error $F_{\lambda,\mu,\nu}(z)$ is smaller when $(\la,\mu,\nu)=(1,1,1)$ than without the algebraic expansion $(\la,\mu,\nu)=(1,1,0)$. For $\theta=0.25\pi$ and $0.50\pi$
(inside the sector $\f{1}{6}\pi<\theta<\f{2}{3}\pi$) it is more accurate with $(\la,\mu,\nu)=(1,1,0)$.
Finally, for $\theta=0.75\pi$ and $0.80\pi$ (inside the sector $\f{2}{3}\pi<\theta<\pi$) only the exponential expansion $E_-(z)$ is present corresponding to $(\la,\mu,\nu)=(0,1,0)$.  It is seen in each case that the first set of $(\lambda,\mu,\nu)$ values, which corresponds to the expansions given in (\ref{e27e}), yields the best approximation.
The results for the remaining values of $\sigma$ can be given a similar interpretation.

Finally, we make a remark about the situation when $\theta$ is such that one of the arguments $ze^{\pm\pi i\sigma}$ of the ${}_1\Psi_0$ functions in (\ref{e51}) is situated on one of the Stokes lines $\pm\pi\kappa$. This corresponds to dealing with some of the borderline cases in (\ref{e27e}) -- (\ref{e29e}). For example, the case $\sigma=\f{1}{6}$ ($\kappa=\f{5}{6}$) and $\theta=\f{2}{3}\pi$, which has $\theta+\pi\sigma=\f{5}{6}\pi$ and $\theta-\pi\sigma=\fs\pi$, corresponds to the borderline case between the second and third expansions in (\ref{e27e}), where $E_+(z)$ is in the process of switching off (as $\theta$ increases). The expansion of ${}_1\Psi_0(xe^{5\pi i/6})$ is given 
by (\ref{e46}).
The expansion of ${}_1\Psi_0(xe^{\pi i/2})$ is given by the first expansion in (\ref{e4exp}), where the algebraic expansion is dominant. Optimal truncation of this last expansion will cancel with the corresponding algebraic contribution from ${}_1\Psi_0(xe^{5\pi i/6})$ in the combination (\ref{e35}). However, the magnitude of the least term in these algebraic expansions is comparable to the exponentially small contributions. Consequently it is not possible within the present framework to deal adequately with such borderline cases. A more refined hyperasymptotic treatment would be necessary to estimate correctly the contribution from the tail of the dominant algebraic expansion.  

A similar remark applies to the Stokes transition of the algebraic expansion that occurs when one of the phases $\theta\pm\pi\sigma=0$. This occurs, for example, when $\sigma=\f{1}{4}$ and $\theta=\f{1}{4}\pi$ corresponding to the borderline case between the first and second expansions in (\ref{e27e}).
However, this case is less acute since this transition of the algebraic expansion is always accompanied by one of the exponential expansions being maximally dominant.

\vspace{0.6cm}

%\begin{center}
{\bf Appendix A: \ An algorithm for the computation of the coefficients $c_j=A_j/A_0$}
%\end{center}
\setcounter{section}{1}
\setcounter{equation}{0}
\renewcommand{\theequation}{\Alph{section}.\arabic{equation}}
We describe an algorithm for the computation of the normalised coefficients $c_j=A_j/A_0$ appearing in the exponential expansion $E_{p,q}(z)$ in (\ref{e22c}). Methods of computing these coefficients by recursion in the case $\alpha_r=\beta_r=1$ have been given by Riney \cite{R} and Wright \cite{W58}; see \cite[Section 2.2.2]{PK} for details. Here we describe an algebraic method for arbitrary $\alpha_r>0$ and $\beta_r>0$. 

The inverse factorial expansion (\ref{e22a}) can be re-written as
\begin{equation}\label{a1}
\frac{g(s)\Gamma(\kappa s+\vartheta')}{\Gamma(1+s)}
=\kappa A_0(h\kappa^\kappa)^{s}\bl\{\sum_{j=0}^{M-1}\frac{c_j}{(\kappa s+\vartheta')_j}+\frac{O(1)}{(\kappa s+\vartheta')_M}\br\}
\end{equation}
for $|s|\to\infty$ uniformly in $|\arg\,s|\leq\pi-\epsilon$, where $g(s)$ is defined in (\ref{e20a}). Introduction of the scaled gamma function $\g^*(z)=\g(z) (2\pi)^{-\fr}e^z z^{\fr-z}$ leads to the representation
\[\g(\alpha s+a)= (2\pi)^\fr e^{-\alpha s} (\alpha s)^{\alpha s+a-\fr} \,{\bf e}(\alpha s; a)\g^*(\alpha s+a),\]
where
\[{\bf e}(\alpha s; a):= e^{-a}\bl(1+\frac{a}{\alpha s}\br)^{\alpha s+a-\fr}=\exp\,\left[(\alpha s+a-\fs) \log\,\left(1+\frac{a}{\alpha s}\right)-a\right].\]

Then, after some routine algebra we find that (\ref{a1}) can be written as
\bee\label{a2}
\frac{g(s) \g(\kappa s+\vartheta')}{\g(1+s)}=\kappa A_0(h\kappa^\kappa)^{s}\,R_{p,q}(s)\,\Upsilon_{p,q}(s),
\ee
where
\[\Upsilon_{p,q}(s):=\frac{\prod_{r=1}^p\g^*(\alpha_rs\!+\!a_r)}{\prod_{r=1}^q\g^*(\beta_rs\!+\!b_r)}\,\frac{\g^*(\kappa s\!+\!\vartheta')}{\g^*(1\!+\!s)},\ \ R_{p,q}(s):=\frac{\prod_{r=1}^p e(\alpha_rs;a_r)}{\prod_{r=1}^q e(\beta_rs;b_r)}\,\frac{e(\kappa s;\vartheta')}{e(s;1)}.\]
Substitution of (\ref{a2}) in (\ref{a1}) then yields the inverse factorial expansion (\ref{e22a}) in the alternative form 
\bee\label{a3}
R_{p,q}(s)\,\Upsilon_{p,q}(s)=\sum_{j=0}^{M-1}\frac{c_j}{(\kappa s+\vartheta')_j}+\frac{O(1)}{(\kappa s+\vartheta')_M}
\ee
as $|s|\to\infty$ in $|\arg\,s|\leq\pi-\epsilon$.

We now expand $R_{p,q}(s)$ and $\Upsilon_{p,q}(s)$ for $s\to+\infty$ making use of the well-known expansion (see, for example, \cite[p.~71]{PK})
\[\g^*(z)\sim\sum_{k=0}^\infty(-)^k\gamma_kz^{-k}\qquad(|z|\ra\infty;\ |\arg\,z|\leq\pi-\epsilon),\]
where $\gamma_k$ are the Stirling coefficients, with 
\[\gamma_0=1,\quad \gamma_1=-\f{1}{12},\quad \gamma_2=\f{1}{288},\quad  \gamma_3=\f{139}{51840},
\quad \gamma_4=-\f{571}{2488320}, \ldots\ .\]
Then we find
\[\g^*(\alpha s+a)=1-\frac{\gamma_1}{\alpha s}+O(s^{-2}),\qquad e(\alpha s;a)=1+\frac{a(a-1)}{2\alpha s}+O(s^{-2}),\]
whence
\[R_{p,q}(s)=1+\frac{{\cal A}}{2s}+O(s^{-2}),\qquad
\Upsilon_{p,q}(s)=1+\frac{{\cal B}}{12s}+O(s^{-2}),\]
where we have defined the quantities ${\cal A}$ and ${\cal B}$ by
\[{\cal A}=\sum_{r=1}^p \frac{a_r(a_r-1)}{\alpha_r}-\sum_{r=1}^q\frac{b_r(b_r-1)}{\beta_r}-\frac{\vartheta}{\kappa} (1-\vartheta),\quad
{\cal B}=\sum_{r=1}^p\frac{1}{\alpha_r}-\sum_{r=1}^q\frac{1}{\beta_r}+\frac{1}{\kappa}-1.\]
Upon equating coefficients of $s^{-1}$ in (\ref{a3}) we then obtain
\bee\label{a4}
c_1=\fs\kappa({\cal A}+\f{1}{6} {\cal B}).
\ee

The higher coefficients are obtained by continuation of this expansion process in inverse powers of $s$. We write the product on the left-hand side of (\ref{a3}) as an expansion in inverse powers of $\kappa s$ in the form
\bee\label{a5}
R_{p,q}(s) \Upsilon_{p,q}(s)=1+\sum_{j=1}^{M-1} \frac{C_j}{(\kappa s)^j}+O(s^{-M})
\ee
as $s\to+\infty$, where the coefficients $C_j$ are determined with the aid of {\it Mathematica}.
From the expansion of the ratio of two gamma functions in \cite[(5.11.13)]{DLMF} we obtain
\[\frac{1}{(\kappa s+\vartheta')_j}=\frac{1}{(\kappa s)^{j}} \bl\{\sum_{j=0}^{M-1}\frac{(-)^k(j)_k}{(\kappa s)^k k!}\,B_k^{(1-j)}(\vartheta')+O(s^{-M})\br\},\]
where $B_k^{(s)}(x)$ are the generalised Bernoulli polynomials defined by
\[\bl(\frac{t}{e^t-1}\br)^{s} e^{xt}=\sum_{k=0}^\infty \frac{B_k^{(s)}(x)}{k!}\,t^k \qquad (|t|<2\pi).\]
Here we have $s=1-j\leq 0$ and $B_0^{(s)}(x)=1$. 

Then the right-hand side of (\ref{a3}) as $s\to+\infty$ becomes
\[1+\sum_{j=1}^{M-1}\frac{c_j}{(\kappa s+\vartheta')_j}+O(s^{-M})=1+\sum_{j=1}^{M-1}\frac{c_j}{(\kappa s)^j} \sum_{k=0}^{M-1} \frac{(-)^k (j)_k}{(\kappa s)^k k!}\,B_k^{(1-j)}(\vartheta')+O(s^{-M})\]
\bee\label{a6}
=1+\sum_{j=1}^{M-1} \frac{D_j}{(\kappa s)^j}+O(s^{-M})
\ee
with
\[D_j=\sum_{k=0}^{j-1} (-)^k\bl(\!\!\!\begin{array}{c}j-1\\k\end{array}\!\!\!\br) c_{j-k}\,B_k^{(k-j+1)}(\vartheta'),\]
where we have made the change in index $j+k\to j$ and used `triangular' summation (see \cite[p.~58]{S}).
Substituting (\ref{a5}) and (\ref{a6}) into (\ref{a3}) and equating the coefficients of like powers of $\kappa s$, we then find $C_j=D_j$ for $1\leq j\leq M-1$,
whence
\[c_j=C_j-\sum_{k=1}^{j-1}(-)^k\bl(\!\!\!\begin{array}{c}j-1\\k\end{array}\!\!\!\br) c_{j-k}\,B_k^{(k-j+1)}(\vartheta').\]

Thus we find 
\begin{eqnarray*}
c_1&=&C_1,\\
c_2&=&C_2-c_1 B_1^{(0)}(\vartheta'),\\
c_3&=&C_3-2c_2 B_1^{(-1)}(\vartheta')+c_1 B_2^{(0)}(\vartheta'), \ldots
\end{eqnarray*}
and so on, from which the coefficients $c_j$ can be obtained recursively.
With the aid of {\it Mathematica} this procedure is found to work well in specific cases when the various parameters have numerical values, where up to a maximum of 100 coefficients have been so calculated.

\vspace{0.6cm}

%\begin{center}
{\bf Appendix B: \ An alternative derivation of the expansion (\ref{e400}) when $\sigma=\fs$}
%\end{center}
\setcounter{section}{2}
\setcounter{equation}{0}
\renewcommand{\theequation}{\Alph{section}.\arabic{equation}}
The associated functions ${}_1\Psi_0(xe^{\pm\pi i\sigma})$  defined in (\ref{e36}) can be expressed in terms of the confluent hypergeometric function ${}_1F_1(a;b;z)$ when $\sigma=\fs$ in the form
\[{}_1\Psi_0(xe^{\pm\fr\pi i})=\sum_{n=0}^\infty \frac{\g(n+\delta)}{(2n)!}\,(-x^2)^n\pm ix\sum_{n=0}^\infty \frac{\g(n+\delta+\fs)}{(2n+1)!}\,(-x^2)^n\]
\bee\label{b1}
=\g(\delta)\,{}_1F_1(\delta;\fs;-\f{1}{4}x^2)\pm ix \g(\delta+\fs)\,{}_1F_1(\delta+\fs;\f{3}{2};-\f{1}{4}x^2).
\ee
The asymptotic expansion of ${}_1F_1(a;b;z)$ on the Stokes lines $\arg\,z=\pm\pi$ has been considered in \cite{PCHF}, where, from Theorem 1 of this reference, it is shown that
\[\frac{\g(a)}{\g(b)}\,{}_1F_1(a;b;-x)=\frac{x^{-a}\g(a)}{\g(b-a)} \sum_{k=0}^{m_o-1} \frac{(a)_k(1+a-b)_k}{k!\,x^k}\]
\bee\label{b2}
+x^\xi e^{-x}\bl\{\cos\,\pi\xi \sum_{j=0}^{M-1}(-)^j c_jx^{-j}-\frac{2\sin\,\pi\xi}{\sqrt{2\pi x}} \sum_{j=0}^{M-1} (-)^jb_jx^{-j}+O(x^{-M})\br\}
\ee
as $x\to+\infty$, where $\xi:=a-b$, $M$ is a positive integer and $m_o\sim x-\Re (2a-b)$ is the optimal truncation index of the algebraic expansion. The coefficients $c_j$ and $b_j$ are given by
\[c_j=\frac{(1-a)_j(b-a)_j}{j!},\qquad b_j=\sum_{k=0}^j (-2)^k (\fs)_k c_{j-k}\,G_{2k,j-k}(\fs),\]
where the first five even-order coefficients $G_{2k,j}(\fs)$ are specified in (\ref{e42c}); higher coefficients are generated from the expansion in (\ref{e49}) with $\mu=1$.
The quantity $\gamma_j$ is given by
\[\gamma_j=m_o+2a-b-x-j=\alpha-j,\]
where $m_o=x-2a+b+\alpha$, with $|\alpha|$ bounded.

We now substitute the expansion (\ref{b2}) into (\ref{b1}) 
%with $\delta=1-b=\vartheta+\fs$. 
and observe that the parameter $\xi=\delta-\fs$
for the first hypergeometric function and $\xi=\delta-1$ for the second, and that the coefficients $c_j$, and hence $b_j$, are the same for both hypergeometric functions. In addition, the optimal truncation index $m_o$ is also the same for the algebraic expansions associated with each of these hypergeometric functions. Then, after some routine algebra, we obtain
\[{}_1\Psi_0(xe^{\pm\fr\pi i})=2e^{\pm\pi i\delta}\sum_{k=0}^{m_o-1} \frac{\g(2k+2\delta)}{k! x^{2k+2\delta}}\hspace{5cm}\]
\[\hspace{3cm}+(Xe^{\pm\pi i})^\vartheta e^{-X}\bl\{\sum_{j=0}^{M-1}\bl\{\frac{1}{2}A_j(\fs)\pm\frac{iB_j(\fs)}{\sqrt{2\pi X}}\br\}(-X)^{-j}+O(X^{-M})\br\}.\]
The variable $X=x^2/4$ and the coefficients
\bee\label{b3}
A_j(\fs)=A_0(\fs) c_j=A_0(\fs)\,\frac{(1-\delta)_j(\fs-\delta)_j}{j!},\qquad B_j(\fs)=A_0(\fs) b_j,
\ee
where, from (\ref{e22b}), $A_0(\fs)=2\sqrt{\pi}$. This agrees with the expansion given in (\ref{e400}) when we put $\delta=1-b$.

%{\bf Received: Month xx, 2017}
\end{document}